\documentclass[12pt]{article}
\usepackage{amssymb,amsfonts,amsmath}
\usepackage[english]{babel}
\usepackage{graphicx}

\begin{document}

\centerline{\large\bf On the convergence of an exotic formal series solution of an ODE}
\bigskip

\centerline{R.\,R.\,Gontsov\footnote{Institute for Information Transmission Problems of RAS,
Bolshoy Karetny per. 19, build.1, Moscow 127051 Russia, gontsovrr@gmail.com.}, 
I.\,V.\,Goryuchkina\footnote{Keldysh Institite of Applied Mathematics of RAS,
Miusskaya sq. 4, Moscow 125047 Russia, igoryuchkina@gmail.com.}}
\bigskip
\bigskip

We consider a non-linear ordinary differential equation (ODE)
\begin{eqnarray}\label{eq1}
F(x,y,\delta y,\ldots,\delta^n y)=0,
\end{eqnarray}
where $\delta=x(d/dx)$ and $F(x,y_0,y_1,\ldots,y_n)$ is a holomorphic function
in  a small polydisc
$$
\Delta=\{|x|<\rho,\; |y_j|<\varepsilon_j,\; j=0,\dots, n\}.
$$
Suppose that the equation \eqref{eq1} possesses a formal solution $y=\varphi,$
\begin{equation}\label{eq2}
\varphi=\sum\limits_{k=0}^{\infty}\alpha_k(x^{{\rm i}\eta})x^k, \qquad {\rm i}=\sqrt{-1}, \quad \eta\in\mathbb{R}\setminus\{0\}, 
\end{equation}
 where $\alpha_k(t)$ are meromorphic functions at the origin:
 $$
\alpha_k(t)=t^{-r_k}\sum\limits_{\ell=0}^{\infty}\alpha_{k\ell}\,t^{\ell}, \qquad t^{r_k}\alpha_k(t)\in\mathbb{C}\{t\}, 
\qquad r_k\in\mathbb{Z},
$$
with some common punctured disc ${\cal D}=\{0<|t|<\varepsilon\}$ of convergence.

The series \eqref{eq2} will be called {\it exotic}, in the terminology of A.\,D.\,Bruno \cite{Br}. In particular, the Painlev\'e III, V, VI equations possess formal solutions of such type \cite{BrGo}, \cite{BrPa}, \cite{Gu}. Thus, our present goal is to obtain some general condition for the convergence of an exotic formal series solution of the equation \eqref{eq1}. In this sense, our work continues a series of articles where similar questions were studied for generalized formal power series \cite{GG1}, \cite{GG3} and formal Dulac series \cite{GG4}, which were inspired by the original paper of B.\,Malgrange \cite{Ma} on the classical formal power series solutions of a non-linear ODE.  
\medskip

{\bf Theorem 1.} {\it Let $(\ref{eq2})$ be a formal solution of the equation $(\ref{eq1}):$
$$
F(x,\Phi)=0, \qquad \Phi:=(\varphi,\delta\varphi,\dots,\delta^n\varphi),
$$
such that $\displaystyle \frac{\partial F}{\partial y_n}(x,\Phi)\not\equiv 0$. Furthermore, let each exotic formal series
$\displaystyle \frac{\partial F}{\partial y_i}(x,\Phi)$ be of the form
$$
\frac{\partial F}{\partial y_i}(x,\Phi)=B_i(x^{{\rm i}\eta})x^N+D_i(x^{{\rm i}\eta})x^{N+1}+\ldots, \qquad i=0,1,\ldots,n,
$$
with the same $N$ for all $i$, $B_n\not\equiv0$ and
$$
{\rm ord}_0\,B_i\geqslant {\rm ord}_0 \,B_n, \qquad i=0,1,\dots,n.
$$
Then the series $(\ref{eq2})$ converges uniformly in some open sector $S\subset{\mathbb C}\setminus{\mathbb R}_+$  with the vertex at the origin and of sufficiently small radius.}
\medskip

We give a sketch of the proof of Theorem 1 in a series of lemmas. The first one is on the reduction of the initial equation to a special form, with the proof in a spirit of Malgrange's proof \cite{Ma} of the corresponding lemma for a classical formal power series solution.
\medskip

{\bf Lemma 1.} {\it Under the assumptions of Theorem 1, there exists $m'\in\mathbb{N}$ such that for any $m\geqslant m'$,  the transformation
\begin{equation}\label{eq3}
 y=\sum\limits_{k=0}^m\alpha_k(t)x^k+x^m u
\end{equation}
reduces the equation \eqref{eq1} to an equation of the form
\begin{equation}\label{eq4}
  x^{{\rm i}\eta\, r}\sum_{i=0}^nA_i(x^{{\rm i}\eta})\,(\delta+m)^iu=x\,M(x^{{\rm i}\eta},x,u,\delta u,\dots,\delta^nu),
\end{equation}
where $r\in{\mathbb Z}_+$, $A_i(t)$ are holomorphic near $t=0$, $A_n(0)\neq 0,\;$ $M(x,t,u_0,\dots,u_n)$ is holomorphic near $0\in{\mathbb C}^{n+3}$.}
\medskip

The reduced equation (\ref{eq4}) possesses an exotic formal series solution
\begin{equation}\label{eq7}
\psi=\sum_{k=1}^{\infty}c_k(x^{{\rm i}\eta\, r})\,x^k, \qquad c_k(t)=\alpha_{k+m}(t)=t^{-\nu_k}\sum_{\ell=0}^{\infty} c_{k\ell}\,t^{\ell},
\end{equation}
moreover, this is its unique formal solution of such a form, since the integer $m$ can be chosen in such a way that for any $k\geqslant1$, a Fuchsian linear differential operator (near $t=0$)
$$
\sum_{i=0}^nA_i(t)\Bigl(k+m+{\rm i}\,\eta\,t\frac{d}{dt}\Bigr)^i
$$
has no integer exponents. To prove Theorem 1, one should prove the convergence of the series $\psi$.

Let us denote 
$$
L(\delta)u:=\sum_{i=0}^nA_i(0)(\delta+m)^iu,
$$
$$
x^{{\rm i}\eta}\, H(x^{{\rm i}\eta},u,\delta u\dots,\delta^nu):=L(\delta)u-\sum_{i=0}^nA_i(x^{{\rm i}\eta})(\delta+m)^iu,
$$
and write down the equation \eqref{eq4} in the form
\begin{equation}
x^{{\rm i}\eta\, r} L(\delta)u=x^{{\rm i}\eta(r+1)} H(x^{{\rm i}\eta},u,\delta u\dots,\delta^nu)+x\,M(x^{{\rm i}\eta},x,u,\delta u,\dots,\delta^nu).
\label{eq10}
\end{equation}

First, by the equation \eqref{eq10} we construct an algebraic equation
\begin{equation}
t^r\,\sigma\,v =t^{r+1} \widetilde{H}(t,v,v,\dots,v)+x\,\widetilde{M}(t,x,v,v,\dots,v), \quad \sigma>0, \label{eq9}
\end{equation}
majorant for the equation \eqref{eq10} in the sense that it has a unique formal solution of the form 
\begin{equation}
\tilde{\psi}=\sum_{k=1}^{\infty}C_k(t)\,x^k, \qquad C_k(t)=t^{-\nu_k}\sum_{\ell=0}^{\infty}C_{k\ell}\,t^{\ell}, \quad C_{k\ell}\geqslant0,  \label{eq11}
\end{equation}
and $|c_{k\ell}|\leqslant C_{k\ell}$ for all $k, \ell$. Let us briefly describe this construction.

The functions $H,M$ of the right hand side in \eqref{eq10} are represented by power series
convergent in $\mathcal{D}\times\Delta$:
\begin{eqnarray*}
H(t,u_0,\dots,u_n)&=&\sum_{i=0}^n\sum_{s=0}^{\infty} h_{s,i}\,t^s u_i, \quad h_{s,i}\in\mathbb C, \\
M(t,x,u_0,\dots,u_n)&=&\sum_{(s,p,Q)\in{\mathbb Z}_+^{n+3}}
\alpha_{s,p,Q}\;t^s\,x^p\,u_0^{q_0}u_1^{q_1} \ldots u_n^{q_n}, \\ & & Q=(q_0,q_1,\ldots,q_n), \quad \alpha_{s,p,Q}\in\mathbb C.
\end{eqnarray*}
Then the functions $\widetilde H, \widetilde M$ in the equation \eqref{eq9} are defined as follows:
\begin{eqnarray*}
\widetilde H(t,v,\dots,v)&=&\sum_{i=0}^n\sum_{s=0}^{\infty}| h_{s,i}|\,t^s v, \\
\widetilde M(t,x,v,\dots,v)&=&\sum_{(s,p,Q)\in{\mathbb Z}_+^{n+3}}
|\alpha_{s,p,Q}|\;t^s\,x^p\,v^{q_0}v^{q_1} \ldots v^{q_n}.
\end{eqnarray*}
The number $\sigma$ is defined by the formula
$$
\sigma=\inf_{{k\in\mathbb N}\atop{\ell\in\mathbb{Z}_+}}|L(k+m+(\ell-kr){\rm i\,}\eta)|.
$$
This value is strictly positive, as $L(k+m+(\ell-kr){\rm i\,}\eta)\ne0\;$ for $\;{k\in\mathbb N},$ $\;{\ell\in\mathbb{Z}_+}$.
\medskip

{\bf Lemma 2.} {\it The equation \eqref{eq9} has a uniquely determined formal solution of the form \eqref{eq11}, which is
majorant for the exotic formal series solution \eqref{eq7} of the equation \eqref{eq10}: $|c_{kl}|\leqslant C_{kl}$.}
\medskip

The proof of this lemma is somehow similar to that of Lemma 2 from \cite{GG2}, where a proof of Malgrange's theorem for a classical power series by the majorant method was proposed.

The second step is to prove the convergence of the series \eqref{eq11}.
Let us write down the equation \eqref{eq9} in the form 
\begin{equation}\label{eq23}
v=\frac{t}{\sigma} \widetilde{H}(t,v,\dots,v)+
\frac{x}{\sigma t^r}\,\widetilde{M}(t,x,v,\dots,v).
\end{equation}
Its right hand side is holomorphic in the domain 
$$
\widetilde{\mathcal{D}}=\{0<\tau<|t|<\varepsilon,\; |x|<\rho,\; |v|<\varepsilon_n\},
$$ 
and we have 
$$
\mathcal{M}=\sup\limits_{\widetilde{\mathcal{D}}}\Bigl|\frac{t}{\sigma} \widetilde{H}(t,v,\dots,v)+
\frac{x}{\sigma t^r}\,\widetilde{M}(t,x,v,\dots,v)\Bigr|<+\infty.
$$

Further we look at the equation \eqref{eq23} as at an equation in two variables $x,$ $v$, with a nonzero parameter $t$. For each fixed value $t=t_0$, $\tau<|t_0|<\varepsilon$, such that
$$
\frac{t_0}{\sigma}\sum_{i=0}^n\sum_{s=0}^{\infty}| h_{s,i}|\,t_0^s\ne1,
$$
we can apply the
implicit function theorem to the equation \eqref{eq23}. According to this theorem, the equation has a unique holomorphic solution $v=v(x,t_0)$ in some open disc 
$$
\mathcal{V}_{t_0}=\Bigl\{|x|<\rho\Bigl(\frac{\varepsilon_n}{\varepsilon_n+2\mathcal{M}_{t_0}}\Bigr)^2\Bigr\},\qquad \mathcal{M}_{t_0}\leqslant\mathcal{M},
$$
whose power series should coincide with the formal series \eqref{eq11}, for $t=t_0$. Hence, the latter converges uniformly in
some open domain
$$
\mathcal{V}=\Bigl\{0<\tau<|t|<\tau',\;|x|<\rho\Bigl(\frac{\varepsilon_n}{\varepsilon_n+2\mathcal{M}}\Bigr)^2\Bigr\},
$$ 
which concludes the proof of Theorem 1.


\begin{thebibliography}{99}

\bibitem{Br}

Bruno, A.\,D., Exotic expansions of solutions to an ordinary differential equation, {\it Doklady Math.} {\bf 76}:2 (2007), 729--733.

\bibitem{BrGo}

Bruno, A.\,D., Goryuchkina, I.\,V., All asymptotic expansions of solutions to the sixth Painlev\'e equation,
{\it Doklady Math.} {\bf 76}:3 (2007), 851--855.

\bibitem{BrPa}

Bruno, A.\,D., Parusnikova, A.\,V., Local expansions of solutions to the fifth Painlev\'e equation, {\it Doklady Math.} {\bf 83}:3 (2011), 348--352.

\bibitem{GG1}

Gontsov, R.\,R., Goryuchkina, I.\,V., On the convergence of generalized power series satisfying an algebraic ODE,
{\it Asympt. Anal.} {\bf 93}:4 (2015), 311--325.

\bibitem{GG2}

Gontsov, R.\,R., Goryuchkina, I.\,V., An analytic proof of the Malgrange theorem on the convergence of formal solutions 
of an ODE, {\it J. Dyn. Control Syst.} {\bf 22}:1 (2016), 91--100.

\bibitem{GG3}

Gontsov, R.\,R., Goryuchkina, I.\,V., The Maillet--Malgrange type theorem for generalized power series, {\it Manuscripta 
Math.} {\bf 156}:1 (2018), 171--185. 

\bibitem{GG4}

Gontsov, R.\,R., Goryuchkina, I.\,V.,  On the convergence of formal Dulac series satisfying an algebraic ODE,
{\it Sbornik Math.}, 2019 (to appear).

\bibitem{Gu}

Guzzetti, D., Tabulation of Painlev\'e 6 transcendents, {\it Nonlinearity} {\bf 25} (2012), 3235--3276.

\bibitem{Ma}

Malgrange, B., Sur le theor\`eme de Maillet, {\it Asympt. Anal.} {\bf 2} (1989), 1--4.







\end{thebibliography}
\end{document}